\documentclass[a4paper,12pt]{article}
\usepackage{comment}
\usepackage{cite}
\usepackage{amsmath,amssymb,amsfonts}
\usepackage[T1]{fontenc}
\usepackage[utf8]{inputenc}
\usepackage{graphicx}
\usepackage{fancyhdr}
\usepackage{float}
\usepackage{xcolor}
\usepackage{authblk}
\usepackage{mathrsfs}
\usepackage{empheq}
\usepackage[hyphens]{url}
\usepackage{hyperref}
\usepackage{amsthm}
\usepackage{tikz}
\usetikzlibrary{decorations.markings}
\usepackage{enumitem}
\usepackage{geometry}

\theoremstyle{plain}

\theoremstyle{definition}

\theoremstyle{remark}

\theoremstyle{conjecture}

\newcommand{\dd}{\mathrm{d}} 

\begin{document}

\title{An Ahmed-like integral}

\author[$\dagger$]{Jean-Christophe {\sc Pain}$^{1,2}$\footnote{jean-christophe.pain@cea.fr}\\
\small
$^1$CEA, DAM, DIF, F-91297 Arpajon, France\\
$^2$Universit\'e Paris-Saclay, CEA, Laboratoire Mati\`ere en Conditions Extrêmes,\\ 
F-91680 Bruy\`eres-le-Châtel, France
}

\date{}

\maketitle

\begin{abstract}
The so-called Ahmed integral,
$$
    \int_{0}^{1}\frac{\arctan\left(\sqrt{2+x^{2}}\right)}{(1+x^{2})\sqrt{2+x^{2}}}\,\dd x=\frac{5\pi^{2}}{96},
$$
has attracted considerable interest since its appearance in the \textit{American Mathematical Monthly} in 2001. Several proofs and extensions have been proposed, including a probabilistic multivariate approach introduced by Pla based on powers of the Gaussian integral. In this note, we extend Pla’s method to the fifth power of the Gaussian integral. By expressing this power as a sequence of iterated integrals and performing successive reductions, we obtain a new integral identity closely related to Ahmed’s integral. In particular, we prove that
$$
    \int_{0}^{1}\frac{\arctan\!\left(\sqrt{\displaystyle\frac{2+x^2}{4+x^2}}\right)}{(1+x^{2})\sqrt{2+x^{2}}}\,\dd x=\frac{\pi^2}{30}.
$$
The derivation suggests that Pla’s technique can systematically generate a family of Ahmed-type integrals associated with higher powers of the Gaussian integral.
\end{abstract}

\section{Introduction}

In 2001, Zafar Ahmed requested the readers of the ``Problems'' section of the \textit{American Mathematical Monthly} to try to evaluate the following rather unusual integral:
\begin{equation*}
    A=\int_{0}^{1}\frac{\arctan\left(\sqrt{2+x^{2}}\right)}{(1+x^{2})\sqrt{2+x^{2}}}\,\dd x.
\end{equation*}
It was shown that $A=5\pi^{2}/96$ \cite{Ahmed2002,oeis}. Since then several authors have provided alternative proofs of this result \cite{Nahin2020}, and even generalisations (see for instance Ref. \cite{Borwein2004}). Ahmed's integral can also be used to calculate an integral due to Coxeter \cite{Coxeter1937,Vidiani2003,Nahin2020}, which reads
$$
    \int_0^{\pi/2}\arccos\left(\frac{\cos x}{1+2\,\cos x}\right)\,\dd x,
$$
and is related to elliptic integrals \cite{Pain2026}. Recently, Pla has proposed an interesting technique to prove Ahmed's integral \cite{Pla2010,Pla2015}. The idea consists in calculating the fourth power of 
\begin{equation*}
    \int_{0}^{+\infty}e^{-x^{2}}\,\dd x=\frac{\sqrt{\pi}}{2},
\end{equation*}
and in expressing it as a multiple (fourth) integral. It turns out that the approach presented in the two latter references can be generalized to higher powers. In this note, we obtain an expression for another Ahmed-like integral, which reads
\begin{equation*}
    B=\int_{0}^{1}\frac{\arctan\left(\sqrt{\displaystyle\frac{2+x^{2}}{4+x^2}}\right)}{(1+x^{2})\sqrt{2+x^{2}}}\,\dd x,
\end{equation*}
and that we found to be equal to 
$$
B=\frac{\pi^2}{30}.
$$

The paper is organized as follows. In Section~2 we briefly review Pla's probabilistic formalism for expressing powers of an integral as iterated multiple integrals, and recall how this method leads to Ahmed's integral when applied to the fourth power of the Gaussian integral. In Section~3 we extend the same approach to the fifth power, which results in a sequence of successive integrations and eventually yields a new integral identity closely related to Ahmed's original formula. A few concluding remarks and possible perspectives for further generalizations are given in the final section.

\section{Pla's probability-integral formalism}

In Refs. \cite{Pla2010,Pla2015}, Pla has explained a calculus method to easily compute the square of a probability integral, and applied it to the computation of Ahmed's integral in Ref. \cite{Pla2015}. Since the technique is general, it seems natural to try to apply it to higher powers. The formalism relies on the following ideas. Let $\alpha$ be a positive number and $f(x,y)$ a function continuous on a subset of $\mathbb{R}^{2}$ containing $(0,\alpha)\times(0,\alpha)$. Then we have
\begin{equation}\label{f1}
    \int_{0}^{\alpha}\int_{0}^{\alpha}f(x,y)\,\dd x\,\dd y=\int_{0}^{1}\,\dd x\int_{0}^{\alpha}\beta\{f(\beta,\beta x)+f(\beta x,\beta)\}\,\dd\beta,
\end{equation}
as well as
\begin{equation}\label{f2}
    \int_{0}^{\alpha}\dots\int_{0}^{\alpha}f(x_{1},\dots,x_{n})\,\dd x_{1}\dots \,\dd x_{n}=\int_{0}^{1}\dots\int_{0}^{1}\,\dd x_{1}\dots \,\dd x_{n-1}\int_{0}^{\alpha}\beta^{n-1}h(\beta, x_1, \cdots, x_n)\,\dd\beta,
\end{equation}
where $h(\beta, x_1, \cdots, x_n)$ is the sum of all the functions $f(\beta x_{1},\dots,\beta x_{p-1},\beta,\beta x_{p+1},\dots,\beta x_{n})$ with $p$ running from 1 to $n$:
$$
    h(\beta, x_1, \cdots, x_n)=\sum_{p=1}^nf(\beta x_{1},\dots,\beta x_{p-1},\beta,\beta x_{p+1},\dots,\beta x_{n}).
$$
Following Pla, let us consider the third power of an integral:
\begin{equation*}
    T(\alpha):=\left(\int_{0}^{\alpha}g(x)\,\dd x\right)^{3}.
\end{equation*}
Taking the derivative of $T$ with respect to $\alpha$, we get easily
\begin{equation*}
    \frac{\dd T}{\dd\alpha}=3g(\alpha)\left(\int_{0}^{\alpha}g(x)\,\dd x\right)^{2}.
\end{equation*}
Making the substitution $x \to \alpha x$ and resorting to formula (\ref{f1}), we deduce that
\begin{equation*}
    \frac{\dd T}{\dd\alpha}=3g(\alpha)\alpha^{2}\left(\int_{0}^{1}g(\alpha x)\,\dd x\right)^{2}=6\int_{0}^{1}\,\dd x\int_{0}^{\alpha}\,\dd\beta\,\alpha^{2}\beta g(\alpha)g(\alpha\beta)g(\alpha\beta x).
\end{equation*}
Then, by integrating with respect to $\alpha$, we get
\begin{equation*}
    \left(\int_{0}^{\alpha}g(x)\,\dd x\right)^{3}=6\int_{0}^{1}\,\dd x\int_{0}^{1}\,\dd\beta\int_{0}^{\alpha}\,\dd\gamma\,\gamma^{2}\beta g(\gamma)g(\gamma\beta)g(\gamma\beta x).
\end{equation*}
By the same technique, but with the help of formula (\ref{f2}) with $n=3$, we can prove that
\begin{equation*}
    \left(\int_{0}^{\alpha}g(x)\,\dd x\right)^{4}=4!\int_{0}^{1}\,\dd x\int_{0}^{1}\,\dd\beta\int_{0}^{1}\,\dd\gamma\int_{0}^{\alpha}\,\dd\delta\,\delta^{3}\gamma^{2}\beta g(\delta)g(\delta\gamma)g(\delta\gamma\beta)g(\delta\gamma\beta x).
\end{equation*}
The latter equation, combined with $g(x)=e^{-x^{2}}$ and taking the limit $\alpha\rightarrow+\infty$, was used by Pla to obtain Ahmed's integral.

\section{Derivation of a new integral}

If we apply the same procedure for the fifth power, we get, \emph{mutatis mutandis}:
\begin{align*}
    \left(\int_{0}^{+\infty}g(x)\,\dd x\right)^{5}=&5!\int_{0}^{1}\,\dd x\int_{0}^{1}\,\dd\beta\int_{0}^{1}\,\dd\gamma\int_{0}^{1}\,\dd\delta\int_{0}^{+\infty}\,\dd\epsilon\,\epsilon^4\delta^{3}\gamma^{2}\beta g(\delta)g(\delta\gamma)g(\delta\gamma\beta)\nonumber\\
    &\qquad\qquad\qquad\qquad\qquad\qquad\qquad\qquad\qquad\qquad\times g(\delta\gamma\beta x)g(\epsilon\delta\gamma\beta x).
\end{align*}
Making the substitution $g(x)=e^{-x^{2}}$, we obtain
\begin{equation*}
    \left(\int_{0}^{+\infty}g(x)\,\dd x\right)^{5}=120\int_{0}^{1}\,\dd x\int_{0}^{1}\,\dd\beta\int_{0}^{1}\,\dd\gamma\int_{0}^{1}\,\dd\delta\int_{0}^{+\infty}\,\dd\epsilon\,\epsilon^4\delta^{3}\gamma^{2}\beta\,e^{-\epsilon^2(1+\delta^2+\delta^2\gamma^2+\delta^2\gamma^2\beta^2+\delta^2\gamma^2\beta^2x^2)}.
\end{equation*}
Noticing that $2\epsilon d\epsilon=d(\epsilon^{2})$, we can integrate with respect to $\epsilon$ and obtain
\begin{equation*}
    \left(\int_{0}^{+\infty}g(x)\,\dd x\right)^{5}=45\int_{0}^{1}\,\dd x\int_{0}^{1}\,\dd\beta\int_{0}^{1}\,\dd\gamma\int_{0}^{1}\,\dd\delta\,\frac{\sqrt{\pi}\beta\gamma^2\delta^3}{(1+(1+(1+(1+x^2)\beta^2)\gamma^2)\delta^2)^{5/2}}.
\end{equation*}
We now integrate with respect to $\delta$ and obtain
\begin{align*}
    \left(\int_{0}^{+\infty}g(x)\,\dd x\right)^{5}=&15\int_{0}^{1}\,\dd x\int_{0}^{1}\,\dd\beta\,\dd\gamma\,\frac{\sqrt{\pi}\beta\gamma^2}{(1+(1+(1+x^2)\beta^2)\gamma^2)^{2}}\nonumber\\
    &\qquad\qquad\qquad\qquad\qquad\times\left[2-\frac{5+3(1+(1+x^2)\beta^2)\gamma^2}{(2+(1+(1+x^2)\beta^2)\gamma^2)^{3/2}}\right].
\end{align*}
For the sake of simplicity, we now integrate with respect to $\beta$ (and not with respect to $\gamma$ as the ``hierarchy'' of variables would have suggested) and obtain
\begin{align*}
    \left(\int_{0}^{+\infty}g(x)\,\dd x\right)^{5}=&15\int_{0}^{1}\,\dd x\int_{0}^{1}\,\dd\gamma\sqrt{\pi}\gamma^2\left[\frac{1-\frac{1}{\sqrt{2+\gamma^2}}}{(1+x^2)\gamma^2(1+\gamma^2)}\right.\nonumber\\
    &\qquad\qquad\qquad\qquad\qquad\qquad\qquad\left.-\frac{1-\frac{1}{\sqrt{2+(2+x^2)\gamma^2}}}{(1+x^2)\gamma^2(1+(2+x^2)\gamma^2)}\right].
\end{align*}
Finally, we integrate with respect to $\gamma$, and find
\begin{align*}
    \left(\int_{0}^{+\infty}g(x)\,\dd x\right)^{5}=&15\frac{\sqrt{\pi}}{12}\int_{0}^{1}\,\dd x\left[\frac{\pi\sqrt{2+x^2}-12\arctan(\sqrt{2+x^2})+12\arctan\left(\sqrt{\displaystyle\frac{2+x^2}{4+x^2}}\right)}{(1+x^2)\sqrt{2+x^2}}\right].
\end{align*}
The right-hand side of the latter expression splits into three parts:
\begin{align*}
    \left(\int_{0}^{+\infty}g(x)\,\dd x\right)^{5}=&5\frac{\sqrt{\pi}}{4}\left[\pi\int_{0}^{1}\frac{1}{1+x^2}\,\dd x-12\int_{0}^{1}\frac{\arctan(\sqrt{2+x^2})}{(1+x^2)\sqrt{2+x^2}}\,\dd x\right.\nonumber\\
    &\qquad\qquad\qquad\qquad\qquad\qquad\qquad\qquad\left.+12\int_{0}^{1}\frac{\arctan\left(\sqrt{\displaystyle\frac{2+x^2}{4+x^2}}\right)}{(1+x^2)\sqrt{2+x^2}}\,\dd x\right].
\end{align*}
The first integral in the right-hand side is equal to

\begin{align*}
    \int_{0}^{1}\frac{1}{1+x^2}\,\dd x=\frac{\pi}{4},
\end{align*}
and the second one is Ahmed's integral:
$$
    \int_{0}^{1}\frac{\arctan\sqrt{2+x^{2}}}{(1+x^{2})\sqrt{2+x^{2}}}\,\dd x=\frac{5\pi^2}{96}.
$$
Taking into account the fact that
$$
    \left(\int_{0}^{+\infty}g(x) \,\dd x\right)^{5}=\frac{\pi^{5/2}}{32},
$$
we obtain
$$
    \frac{\pi^{5/2}}{32}=5\frac{\sqrt{\pi}}{4}\left[\frac{\pi^2}{4}-12\times\frac{5\pi^2}{96}+12\int_{0}^{1}\frac{\arctan\left(\sqrt{\displaystyle\frac{2+x^2}{4+x^2}}\right)}{(1+x^2)\sqrt{2+x^2}}\,\dd x\right].
$$
or equivalently
$$
    \frac{\pi^{2}}{8}=5\left[\frac{\pi^2}{4}-\frac{5\pi^2}{8}+12\int_{0}^{1}\frac{\arctan\left(\sqrt{\displaystyle\frac{2+x^2}{4+x^2}}\right)}{(1+x^2)\sqrt{2+x^2}}\,\dd x\right],
$$
which gives
$$
    \frac{\pi^{2}}{8}=5\left[-\frac{3\pi^2}{8}+12\int_{0}^{1}\frac{\arctan\left(\sqrt{\displaystyle\frac{2+x^2}{4+x^2}}\right)}{(1+x^2)\sqrt{2+x^2}}\,\dd x\right].
$$
One thus gets simply
$$
    12\int_{0}^{1}\frac{\arctan\left(\sqrt{\displaystyle\frac{2+x^2}{4+x^2}}\right)}{(1+x^2)\sqrt{2+x^2}}\,\dd x=\frac{\pi^{2}}{40}+\frac{3\pi^2}{8}=\frac{16\pi^2}{40}=\frac{2\pi^2}{5},
$$
yielding the final result
$$
    \int_{0}^{1}\frac{\arctan\left(\sqrt{\displaystyle\frac{2+x^2}{4+x^2}}\right)}{(1+x^2)\sqrt{2+x^2}}\,\dd x=\frac{\pi^2}{30}.
$$

\section{Conclusion}

We have obtained, following a probabilistic multi-variate technique published by Pla a few years ago, and applied by the author to Ahmed's integral, an integral which, to our knowledge, was never obtained before. To derive that result, we just pushed forward Pla's method to the next order, which corresponds to the fifth power of the Gauss integral, and thus to $5-1=4$ successive integrations. Although the calculations become rather cumbersome, they remain feasible, which suggests that a general class of integrals could be obtained, relying on $n-1$ successive integrations, and the fact that
$$
\left(\int_{0}^{\alpha}g(x)\,\dd x\right)^{n}=\frac{\pi^{n/2}}{2^n},
$$
If the observed patterns remain, it might be that only the argument of the $\arctan$ function will be modified, but this has to be checked. Moreover, the fact that Ahmed's integral appears in the course of our derivation of the new integral seems to indicate that an inductive method would be suitable. This will be the subject of a future work.


\begin{thebibliography}{99}

\bibitem{Ahmed2002}
Z. Ahmed, Definitely an Integral, {\it Am. Math. Mon.} {\bf 109}, 670-671 (2002).

\bibitem{Borwein2004} 
J. M. Borwein, D. H. Bailey and R. Girgensohn, \textit{Experimentation in Mathematics: Computational Paths to Discovery}, AK Peters Ltd (March 2004), pp. 17--20.

\bibitem{Coxeter1937}
H.~S.~M.~Coxeter,
\newblock Some integrals connected with trigonometry,
\newblock {\em Math. Gaz.}, \textbf{21}, 55--57 (1937).

\bibitem{Nahin2020}
P. J. Nahin, {\it Inside interesting integrals}, Springer Nature Switzerland AG, 2020.

\bibitem{oeis} 
OEIS Foundation Inc. (2026), Entry A005132 in The On-Line Encyclopedia of Integer Sequences, \url{https://oeis.org/A096615}

\bibitem{Pain2026}
J.-C. Pain, {\it Elliptic integral identities derived from Coxeter's integrals}, arXiv:2603.04637, \url{https://arxiv.org/abs/2603.04637} (2026).

\bibitem{Pla2010} 
J. Pla, A footnote to the theory of double integrals, {\it Math. Gaz.} {\bf 94}, 262--269 (2021).

\bibitem{Pla2015}
J. Pla, {\it A tale of two integrals: The probability and Ahmed's integrals}, arXiv:1505.03314, \url{https://arxiv.org/abs/1505.03314} (2015).

\bibitem{Vidiani2003}
G.~Vidiani,
\newblock Les int\'egrales de Coxeter,
\newblock {\em Quadrature, Magazine de math\'ematiques pures et \'epic\'ees}, no.~50,
Oct.–Dec.~2003, pp.~7--12 (in French).

\end{thebibliography}
\end{document}